\documentclass[a4paper,twoside,final]{siamltex}

\usepackage[a4paper,top=3cm,bottom=2.5cm,left=3cm,right=2.5cm]{geometry}

\usepackage{amsmath}
\usepackage{etex}
\usepackage{amsopn}
\usepackage{float}
\usepackage{graphicx}
\usepackage{booktabs}
\usepackage{amsfonts}

\usepackage{algorithm}
\usepackage{algpseudocode}
\usepackage{xpatch}
\usepackage[percent]{overpic}
\usepackage{tikz}
\usetikzlibrary{spy}
\usepackage{bm}

\usepackage{xcolor}
\usepackage[british]{babel}
\usepackage{lmodern}
\usepackage{multirow}
\usepackage{mwe}
\usepackage{siunitx}
\usepackage{textcomp}
\usepackage[unicode=true,
 bookmarks=true,bookmarksnumbered=false,bookmarksopen=false,
 breaklinks=false,pdfborder={0 0 0},backref=false,colorlinks=false]
 {hyperref}
\usepackage[nameinlink]{cleveref}
\crefname{figure}{figure}{figures}
\newlength{\tempdima}
\newcommand{\rowname}[1]
{\rotatebox{90}{\makebox[\tempdima][c]{\textbf{#1}}}}

\algrenewcommand\alglinenumber[1]{{\sffamily\footnotesize#1}}

\makeatletter
\xpatchcmd{\algorithmic}{\itemsep\z@}{\itemsep=0.5ex plus1pt}{}{}
\makeatother

\newcommand{\Rbb}{\mathbb{R}}

\DeclareMathOperator*{\argmin}{arg\,min}
\DeclareMathOperator{\TV}{TV}
\DeclareMathOperator{\dTV}{dTV}
\DeclareMathOperator{\prox}{prox}
\DeclareMathOperator{\tol}{tol}

\newcommand{\func}[3]{%
	#1\colon  #2  \to  #3
}

\usepackage{doi}

\def\PicWidth{1cm}
\definecolor{salmon}{rgb}{1.0, 0.55, 0.41}
\definecolor{frames}{rgb}{0.13, 0.7, 0.67}

\newcommand{\bx}{\bm{u}}
\newcommand{\bb}{\bm{b}}
\newcommand{\bv}{\bm{v}}

\newcommand{\bu}{\bm{u}}
\newcommand{\by}{\bm{y}}
\newcommand{\bz}{\bm{z}}
\newcommand{\bq}{\bm{q}}

\algrenewcommand\alglinenumber[1]{{\sffamily\footnotesize#1}}

\makeatletter
\xpatchcmd{\algorithmic}{\itemsep\z@}{\itemsep=0.5ex plus1pt}{}{}
\makeatother

\usepackage{xargs}                      
\begin{document}

\title{Synergistic Multi-spectral CT Reconstruction with Directional Total Variation}

\author{Evelyn Cueva\footnotemark[1] \and Alexander Meaney\footnotemark[2] \and Samuli Siltanen\footnotemark[2] \and Matthias J. Ehrhardt\footnotemark[3]}

\renewcommand{\thefootnote}{\fnsymbol{footnote}}

\footnotetext[1]{Research Center on Mathematical Modelling (MODEMAT), Escuela Polit\'ecnica Nacional, Quito, Ecuador (ecueva@uchile.cl)}

\footnotetext[2]{Department of Mathematics and Statistics, University of Helsinki,
Finland}
\footnotetext[3]{Institute for Mathematical Innovation, University of Bath,
United Kingdom}

\renewcommand{\thefootnote}{\arabic{footnote}}

\pagestyle{plain} \thispagestyle{plain} \markboth{
E.~CUEVA} {Synergistic Multi-spectral CT Reconstruction with Directional Total Variation}

\maketitle

\begin{abstract}
  This work considers synergistic multi-spectral CT reconstruction where information from all available energy channels is combined to improve the reconstruction of each individual channel, we propose to fuse this available data (represented by a single sinogram) to obtain a polyenergetic image which keeps structural information shared by the energy channels with increased signal-to-noise-ratio. This new image is used as prior information during a channel-by-channel minimization process through the directional total variation. We analyze the use of directional total variation within variational regularization and iterative regularization. Our numerical results on simulated and experimental data show improvements in terms of image quality and in computational speed.

\end{abstract}


\maketitle

\section{Introduction}
\subsection{Undersampling in multi-spectral CT}
Computed tomography (CT) is a widely used technique in many different fields of science and industry; for example in medicine, it enables visualizing the internal structure of a patient. The principle of this technique is to study the attenuation of X-rays when they pass through the target object~\cite{natterer2001mathematics, de2014industrial}. Despite the potential usefulness of CT, the X-ray source produces a single energy spectrum and the detector does not discriminate between photon energies. As a consequence, two tissues whose elemental composition are different might be indistinguishable in the resulting CT image~\cite{shikhaliev2008energy, shikhaliev2011photon,mccollough2015dual}. The latter makes it difficult to identify and classify different tissues and motivates the multi-spectral techniques based on new scanner technologies~\cite{macovski1976energy, alessio2013quantitative,marin2014state}.

Dual and multi-spectral CT use different technical approaches for acquiring multi-energetic data, \emph{e.g.}, the rapid tube potential switching, multilayer detectors, or dual (multi) X-ray sources~\cite{zou2008analysis,xu2009dual,krauss2011dual, johnson2012dual}. This multi-energetic data provides much more information about the tissue composition allowing to differentiate its constituent materials~\cite{mccollough2015dual} but, in addition, the measurement process needs a balance between radiation dose, acquisition time and image quality. A reduction in radiation dose is achieved by reducing the number of views in the acquisition which, in turn, decreases the spatial resolution~\cite{mccollough2012achieving, hamalainen2013sparse,lell2015evolution}. A recent study proposes to reconstruct a multi-spectral CT image by reducing the dose in each energy window, when just a limited and non-overlapping range of angles is observed~\cite{toivanen2019joint}. As the resolution of the reconstructions is affected by this lack of measurements, small objects
cannot be reconstructed and the resulting images are affected by the presence of artifacts~\cite{huesman1977effects, frikel2013characterization}. Many advanced reconstruction techniques have been proposed to simultaneously or independently reconstruct an spectral-image in this scenario~\cite{gao2011multi,kazantsev2018joint,niu2018nonlocal,
wu2018spatial,toivanen2019joint,hu2019nonlinear}. For example, variational methods are commonly used since they allow to directly incorporate prior information and constraints into the model~\cite{scherzer2009variational}. In addition, regularizers can be added as part of the objective functional or in the optimization process, to overcome ill-posedness~\cite{sarkar1981some}. The expected structural correlation between different energy levels has motivated the use of structural priors to improve joint or individual reconstructions~\cite{kaipio1999inverse, rigie2015joint, kazantsev2018joint, hintermuller2018function, toivanen2019joint}, some of them based on level sets methods~\cite{ehrhardt2013vector}. For example, in~\cite{kazantsev2018joint}, the authors proposed a joint image reconstruction method for x-ray spectral computed tomography. Their proposed dTV-p method propagates the structural information along channels using direction total variation (dTV) with a reference image chosen by a probability mass function and considering a wide number of channels (70) for a synthetic phantom. Direction total variation has been successfully used in several other medical imaging applications~\cite{ehrhardt2014joint, ehrhardt2016multicontrast, ehrhardt2016pet} as in hyperspectral remote sensing~\cite{bungert2018blind}. Another similar regularization applied in joint reconstructions to include structural dependence between channels, is the total nuclear variation (TVN) presented in~\cite{rigie2015joint, rigie2017assessment}.


\subsection{Main contribution}

We propose a novel reconstruction technique to solve the undersampling problem in multi-spectral CT, where information from all available energy channels is combined to obtain a polychromatic image. The latter keeps the structural information shared by the energy channels and is used to improve the reconstruction of each individual channel (channel-by-channel reconstruction) using directional total variation (dTV). We explore variational and iterative regularization methods, specifically, the forward-backward splitting algorithm (FBS)~\cite{combettes2005signal,combettes2011proximal} and Linearized Bregman iterations~\cite{osher2005iterative,yin2008bregman,yin2010analysis, benning2017choose,corona2019enhancing} to solved the undersampling problem using simulated and experimental data. The combination of these methods and dTV shows improvements in terms of image quality and computational speed compared, for example, to joint reconstruction techniques as TNV.

In section 2, we describe the inverse problem behind multi-spectral CT data seen as a minimization problem. Later, in section 3, we present the variation and iterative regularization of the inverse problem and we describe the FBS algorithm and Bregman iterations to solved them,respectively. We include total variation and directional total variation regularizers needed during the regularization process. The last section is devoted to present numerical results using synthetic and real data. Here, we specify all the settings needed during the reconstruction. 

\section{Inverse problem}

Multi-spectral CT aims to recover energy-dependent attenuation maps $\bx_k$ of a target object for energies $E_k$ with $k=1,...,K$.
The acquisition method considers X-ray projections using only a limited set of angles, \emph{i.e.}, we want to reconstruct $\bx_k \in  \Rbb^N$ given data $\bm b_k\in \Rbb^{M}$ where $M \ll N$. When a considerable amount of 
measurements is available, classical methods such as filtered back projection, Kaczmarz iterations or iterative techniques can be used to solve an associated linear system of the form $A\bx_k = \bb_k$ or the associated least squares problem (see \emph{e.g.}~\cite{ natterer2001mathematics,clackdoyle2010tomographic})
\begin{equation}\label{eq:data_fit_minimization_1v}
\min_{\bx_k\in\Rbb^N}\frac 12 \|A\bx_k -\bb_k\|^2_2,
\end{equation}
where $A$ is the forward operator (a matrix in the discrete case) that relates the image $\bu_k$ to the given data $\bb_k$. The ill-posedness of this inverse problem makes a direct inversion of the matrix $A$ unstable even for a suitable number of measurements. The undersampling
scenario is even more challenging, since $M\ll N$, the system is under-determined. 

For 2D CT, $M = m_{1}\cdot m_{2}$
where $m_{1}$ is the number of angles and $m_{2}$ is the number of detectors, and $N = n_1\cdot n_2$, where $n_1$ and $n_2$ are the number of rows and
columns of $\bx_k$ (considered as a matrix of pixels), respectively. 

\subsection{Forward model}

We recall the forward modelling for multi-spectral CT. For a fixed energy channel $E_k$, an initial intensity $I_i^0(E_k)$ of X-rays is emitted along a line $L_i$ (from source to detector) given a final intensity $I_i^1$, for $i=1,\ldots M$. The discretized linear model use for reconstructing a vectorized image $\bu(E_k)$ of $N$ pixels (see, \emph{e.g}~\cite{toivanen2019joint}) is given by
\begin{equation}\label{eq:measurements}
  b_{ik} := -\ln\left(\frac{I^1_i}{I^0_i(E_k)}\right) \approx \sum_{j=1}^N a_{ij}u_j(E_k).
\end{equation}
In~\eqref{eq:measurements}, $u_j(E_k)$ is the value of $\bu(E_k)$ in the corresponding pixel $j$, and $a_{ij}$ is the length of the intersection of the $i$-th line and the $j$-th pixel.

Based on the discretization presented above for each energy $E_k$, we establish the forward model for the projection data 
\begin{equation}\label{eq:linear_system_k}
\bm{b}_k = A_k \bx_k + \bm{e}_k  ,\qquad k =1,\ldots, K, 
\end{equation}
where $A_k\in \Rbb^{M\times N}$, is a matrix with components $a_{ij}$, the vector $\bx_k\in \Rbb^{N}$ has components $u_j$, $\bm{e}_k$ models the measurements noise and, $\bm{b}_k$ is the vector of measurements $b_{ik}$ in \eqref{eq:measurements} for the fixed energy level $k$. The matrix $A_k$ represents the discretization of
the X-ray transform for a particular projection geometry.

From now on, we omit the energy sub-index in~\eqref{eq:linear_system_k} since we will solve an independent problem for each energy channel. 


\section{Regularization}

In this section we discuss the regularizers used in this work, total variation and directional total variation, and how these can be used to regularize an inverse problem. To this end we consider variational regularization and iterative regularization based on Bregman iterations.

\subsection{Regularizers}\label{sec:regularizers}

\subsubsection{Total Variation}
The total variation (TV) regularization has been widely studied due to its 
edge-preserving properties~\cite{rudin1992nonlinear}. TV is well-known to promote piecewise constant images with sharp edges. 

To define the (discrete) total variation for an image $\bm u$ of $n_1\times n_2$ pixels, \emph{i.e.}, $\bm u\in \Rbb^{N}$ with $N=n_1\cdot n_2$, we first introduced the discrete gradient $\func \nabla {\Rbb^{N}}{(\Rbb^2)^N}$ based on a finite difference scheme acting on the image pixels, as follows
\[
\left((\nabla\bm u)_{j}\right)_1=\left\{
\begin{array}{ll}
    u_{j+r(j)}- u_{j}, & \text{if } r(j) <N   \\
    0, & \text{otherwise}
\end{array}
\right.\qquad 
\left((\nabla\bm u)_{j}\right)_2=\left\{
\begin{array}{ll}
    u_{j + s(j)} - u_{j}, & \text{if }  s(j)<N   \\
    0, & \text{otherwise}
\end{array}
\right.
\]
where $r(j)$ is the right neighbour of pixel $j$, and $s(j)$ is the neighbour below pixel $j$. Then
\[
\TV(\bm u) = \sum_{j=1}^{N}\left(\left[\left((\nabla\bm u)_{j}\right)_1\right]^2+ \left[\left((\nabla\bm u)_{j}\right)_2\right]^2 \right)^{1/2}.
\]

\subsubsection{Directional Total Variation}
While TV is a powerful regularizer, it is unclear how additional structural a-priori information can be included. To this end we utilize the directional total variation (dTV) proposed in~\cite{ehrhardt2016multicontrast}.
Let $\bm{\xi}\in (\Rbb^2)^N$ be a vector field with $\|\bm{\xi}_i\|\leq \eta < 1$. We denote by $\mathbf{P}\in \left(\Rbb^{2\times 2}\right)^{N}$, $\mathbf{P}_i := \mathbf{I}-\bm{\xi}_i\otimes\bm{\xi}_i$ an associated matrix-field, where $\mathbf{I}$ is the $2\times 2$ matrix and $\otimes$ represents 
the outer product of vectors. Then $\func \dTV {\Rbb^N}{\Rbb}$ is defined as
\begin{equation}\label{eq:dTV}
\dTV(\bx; \bv) = \sum_j\|\mathbf{P}_j(\nabla \bx)_j\|,
\end{equation}
where $\mathbf{P}_j$ implicitly depends on $\bv$ by means of $\bm{\xi}$.

Some interpretations of dTV are detailed in \cite{ehrhardt2014joint,bungert2018blind}.
We briefly describe some useful properties of this functional using the explicit
expression $\mathbf{P}_j(\nabla \bx)_j = (\nabla \bx)_j -\langle \bm{\xi}_i, (\nabla \bx)_j\rangle\bm{\xi}_i$.
We observe two particular cases:
\[
\mathbf{P}_j(\nabla \bx)_j =
\begin{cases}
(1-\|\bm{\xi}_j\|^2)(\nabla \bx)_j, & \text{if } (\nabla \bx)_j \text{ is parallel to } \bm{\xi}_j \\ 
(\nabla \bx)_j, &  \text{if } (\nabla \bx)_j \text{ is perpendicular to } \bm{\xi}_j.
\end{cases}
\]
So, when we minimize $\dTV(\bx)$, we are favouring $\bx$ such that its gradient is collinear to the direction $\bm{\xi}_i$ as long as $\|\bm{\xi}_i\|\not =0$. We note that a vanishing gradient $\nabla u = 0$ always leads to a smaller function value such that no artificial jumps are enforced. 

In order to incorporate dTV into our model, we define the vector field below based
on the known image $\bv\in \Rbb^N$ by
\begin{equation}\label{eq:eta}
  \bm{\xi}_j =\eta \frac{(\nabla \bv)_j}{\|(\nabla \bv)_j\|_\varepsilon}
\end{equation}
with $\|\bx\|_\varepsilon = \sqrt{\|\bx\|^2+\varepsilon^2}$. The parameter $\varepsilon>0$ avoids singularities when $(\nabla \bv)_j =0$ and $\eta$ is an edge parameter related to the size of an edge. 

\begin{figure}[htbp]
\centering\begin{tabular}{@{}c@{ }c@{ }c@{ }c@{}}
  image & side information & TV & dTV \\
\includegraphics[width=0.25\linewidth]{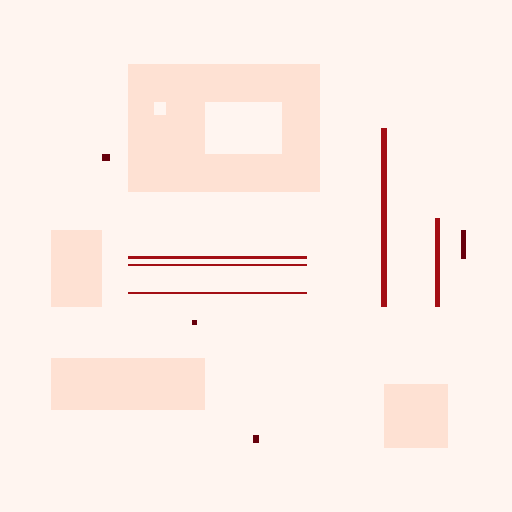}&
\includegraphics[width=0.25\linewidth]{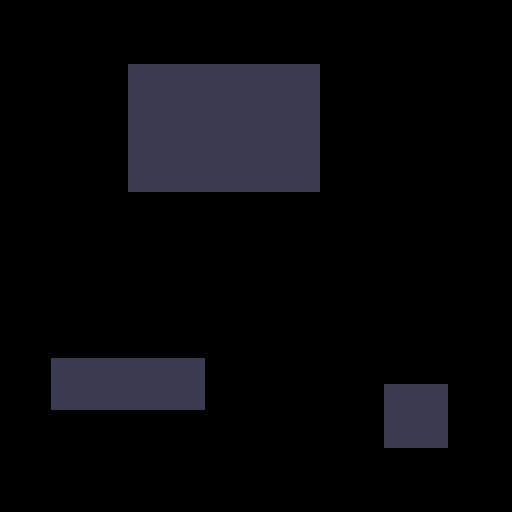}&
\includegraphics[width=0.25\linewidth]{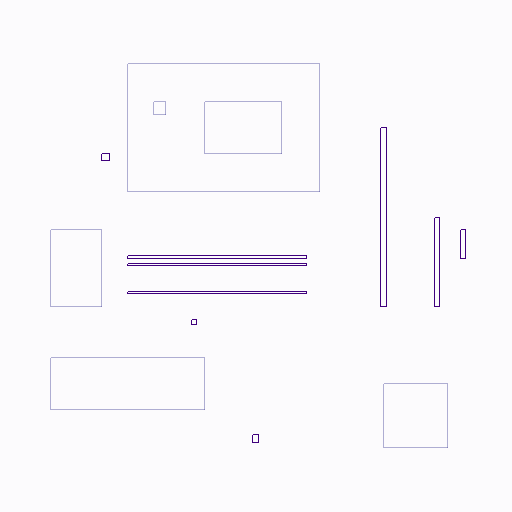}&
\includegraphics[width=0.25\linewidth]{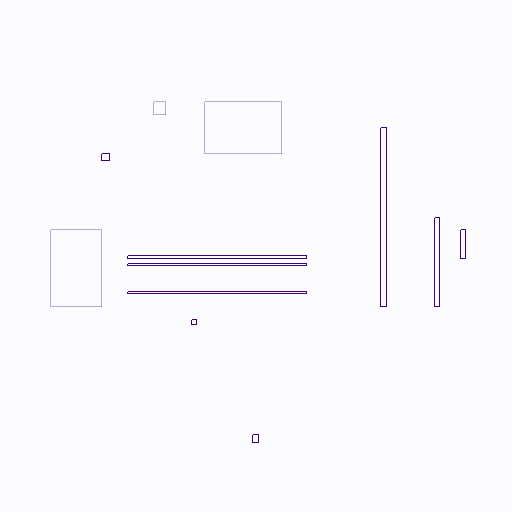}\\
\end{tabular}
\caption{\emph{From left to right:} An image $\bu$, side information $\bv$, pointwise TV-norm $j\mapsto \|(\nabla u)_j\|$, and pointwise dTV-norm $j\mapsto \|\mathbf{P}_j(\nabla \bx)_j\|$ as in~\eqref{eq:dTV}.}\label{fig:TVvsdTV}
\end{figure}

Figure~\ref{fig:TVvsdTV} shows an example which compared TV and dTV. In contrast to TV, dTV only penalizes edges which are missing in the side information.

\subsection{Variational regularization} \label{sec:variational_methods}

A strategy to reconstruct $\bx := \bx_k$ in~\eqref{eq:linear_system_k} is to solve
\begin{equation}\label{eq:variational_1v}
\bx^\ast \in  \argmin_{\bx\in \Rbb^N} \left\{\frac 12\|A\bx-\bb\|^2_2 +\alpha J(\bx) + \iota_{[0,\infty)^N}(\bx)\right\}.
\end{equation}
The first term in~\eqref{eq:variational_1v} is called the data-fit and forces $A\bx$ to stay close to the data, and the regularizer $J$ promotes stability of the inversion. The parameter $\alpha>0$ balances the data-fit term and the regularization provided by $J$. We
will use TV and dTV as $J$. The additional term $\iota_{[0,\infty)^N}(\bx)$ is included to impose a nonnegative constraint for each component of the solution $\bx^\ast$ and is defined as:
\[
\iota_{[0,\infty)^N}(\bx) =
\begin{cases}
0, & \text{ if } u_j\geq 0 \\ 
\infty, & \text{otherwise.}
\end{cases}
\]
Depending on the type of regularization that we choose, we define the following functions:
\begin{align}
\label{eq:G_TV} G_{\TV}(\bx)  &{}= \alpha \TV(\bx) + \iota_{[0,\infty)}(\bx),\\
\label{eq:G_dTV} G_{\dTV}(\bx) &{}= \alpha \dTV(\bx,\bv) + \iota_{[0,\infty)}(\bx).
\end{align}
\subsubsection{Forward-backward splitting algorithm}
The forward-backward splitting (FBS) algorithm~\cite{combettes2005signal} solves the composite minimization problem
\begin{equation}\label{eq:fbs_min_function}
\min_{\bx} \left\{F(\bx) + G(\bx)\right\}
\end{equation}
where $\func F X \Rbb$ and $\func G X {(-\infty, \infty]}$ are two proper, lower semi-continuous and convex functionals such that $F$ is differentiable on $X$ with a $L$-Lipschitz continuous gradient for some $L \in (0,\infty)$.

The principle of this algorithm is based on the two following steps:
\begin{enumerate}
    \item\label{item:fbs_i} a forward (explicit) gradient step on $F$, \emph{i.e.}  $\bx^{t+1/2} = \bx^t - \sigma^t \nabla F(\bx^t)$, and
    \item\label{item:fbs_ii} a backward (implicit) step involving only $G$, \emph{i.e.}  $\bx^{t+1} = \prox_{\sigma^t G} (\bx^{t+1/2}),$
     where the proximal operator is given by
    \begin{equation}\label{eq:proximal_minimization}
    \prox_{\sigma G}(\bz) = \argmin_{\by} \left\{ \frac 12 \|\by-\bz\|^2 +  \sigma G(\by)\right\}.
    \end{equation}
\end{enumerate}
The step size $\sigma^t$ is chosen in each iteration so that it satisfies the
descent inequality
\begin{equation}\label{eq:descent_inequality}
F\left(\bx^{t+1}\right) \leq F\left(\bx^t\right) + \langle \nabla F\left(\bx^t\right), \bx^{t+1} - \bx^t \rangle + \frac{1}{2\sigma^t}\|\bx^{t+1}-\bx^t\|^2.
\end{equation}
More precisely, we reduce $\sigma^t$ until the condition~\eqref{eq:descent_inequality} is satisfied. This
selection of $\sigma$ is known as \emph{backtracking} and is considered in FBS and Bregman iterations.

Now, comparing problem~\eqref{eq:fbs_min_function} with \eqref{eq:variational_1v}, we choose the functions $F$ and $G$ as
\[
  F(\bu) = \frac 12\|A\bu-\bb\|^2_2,\qquad G(\bu)= \alpha J(\bu) + \iota_{[0,\infty)^N}(\bu).
\]
We use the Fast Gradient Projection (FGP) algorithm presented in~\cite{beck2009fast}, to compute the proximal operator of the regularization functionals TV and dTV. For this latter, the detailed algorithm was presented in~\cite[Algorithm 1]{ehrhardt2016multicontrast} where dTV was introduced for the first time. As described in~\cite{ehrhardt2016multicontrast}, the minimization problem associated to the definition of proximal operator is dualized, this dual variable $p$ is initialized as zero and following Algorithm 1 in~\cite{ehrhardt2016multicontrast}, we set the number of (inner) iterations to be 200 and a tolerance of $10^{-5}$, the first condition that is reached stops the algorithm. 
Additionally, we define the objective function value at point $\bu$ as
$H(\bu) = F(\bu) + G(\bu).$ Since FBS algorithm converges to a minimizer of $H$~\cite{combettes2005signal}, we stop the algorithm when the
difference between two consecutive iterations of the $H$ value is less than a given tolerance $\tol$, \emph{i.e.} $H(\bu^{t+1})-H(\bu^{t})\leq \tol\cdot H(\bu^{t+1})$. The~\cref{alg:fbs} describes one iteration of FBS algorithm.

\begin{algorithm}[htbp]
\begin{algorithmic}[1]
\State $\bu^{t+1} = \prox_{\sigma^t G}\left(\bu^t - \sigma^t\nabla F(\bu^t)
\right).$
\If {$F(\bu^{t+1}) > F(\bu^t) + \langle \nabla F(\bu^t), \bu^{t+1} - \bu^t \rangle + \frac{1}{2\sigma^t}\|\bu^{t+1} - \bu^t\|^2$}
\State $\sigma^t = \underline \rho \sigma^t$, for any $\underline \rho <1$ and go back to Step 2.
\Else \State $\sigma^{t+1} = \overline \rho\sigma^t$, for any $\overline \rho >1$.
\EndIf
\end{algorithmic}
\caption{An iteration of forward-backward splitting algorithm}\label{alg:fbs}
\end{algorithm}

\subsection{Iterative regularization}
\begin{figure}[htbp]
  \includegraphics[width=\linewidth]{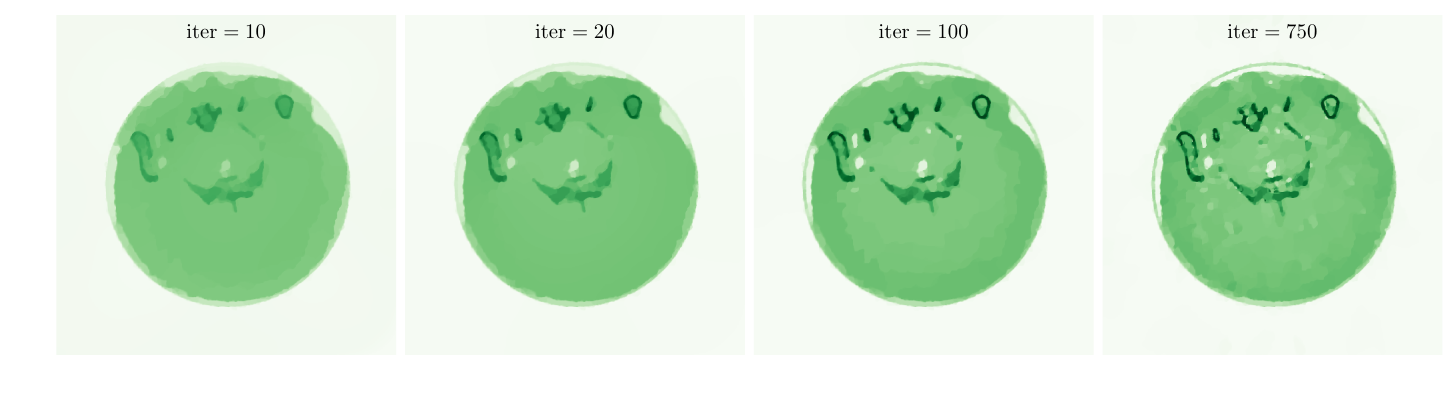}
  \caption{Iterations along the Linearized Bregman iterations. While early iterations are very smooth, the iterates become gradually better defined and eventually the measurement noise is introduced.}\label{fig:bird_semiconvergence}
\end{figure}

A different way to achieve regularization is to apply an iterative method to directly solve
the problem~\eqref{eq:data_fit_minimization_1v}. Iterative methods start with a some vector $\bx^0$ and generate a sequence $\bx^1,\bx^2,\ldots$ that converges to some solution. Usually in these methods, initial iterates $\bx^t$ are fairly close to the exact solution. However, for later iterations, the solutions start to diverge from the desired one and tend to converge to the naive solution $A^{-1}\bb$. Thus, the success of these methods relies on stopping the iterations at the right time. This behavior is known as \emph{semiconvergence}~\cite{hansen2010discrete,natterer2001mathematics} and it is a frequently used tool to solve large-scale problems. In~\cref{fig:bird_semiconvergence} we present an example of this effect. Additionally, iterative regularization avoids a predetermined regularization parameter, and instead, the number of iterations takes the role of a regularization parameter~\cite{hansen2010discrete}. This is an advantage compared to variational regularization since in this latter, a minimization problem needs to be solved every time that a new regularization parameter $\alpha$ is tested. 

\subsubsection{Linearized Bregman iterations}\label{sec:bregman_section}
Under this group of iterative methods, we explore the Linearized Bregman iterations. This algorithm allows us to solve the least squares problem~\eqref{eq:data_fit_minimization_1v}.

We consider the (linearized) Bregman iterations~\cite{osher2005iterative,  yin2010analysis} which makes use of the Bregman distance defined in terms of a given functional $J$ by
\[
D^{\bq^t}_J(\bu,\bu^t)=J(\bu)-J(\bu^t)-\langle \bq^t, \bu-\bu^t\rangle,
\]
where $\bq^t \in \partial J(\bu^t)$ is an element of the sub-differential of $J$ at point $\bu^t$. 

Linearized Bregman iterations are defined as
\begin{align}
\label{eq:bregman_1v_u}\bu^{t+1} &{}= \prox_{\sigma^tG}\left(\bu^t + \sigma^t(\bq^t-\nabla F(\bu^t))\right)\\
\label{eq:bregman_1v_q}\bq^{t+1} &{}=  \bq^t-\frac 1{\sigma^t}\left(\bu^{t+1}-\bu^t+\sigma^t\nabla F(\bu^t)\right)
\end{align}
where $F(\bu) = \frac 12 \|A\bu-\bb\|^2_2$ is the objective function value and $G$ can be chosen as $G_{\TV}$ or $G_{\dTV}$ from~\eqref{eq:G_TV} and \eqref{eq:G_dTV}, respectively. The algorithm with backtracking is detailed in~\cref{alg:bregman_1v}.

\begin{algorithm}[htbp]
\begin{algorithmic}[1]
  \State $\bu^{t+1} = \prox_{\sigma^t G}(\bu^t + \sigma^t(\bq^t-\nabla F(\bu^t)))$
  \If {$F(\bu^{t+1}) > F(\bu^t) + \langle \nabla F(\bu^t), \bu^{t+1} - \bu^t \rangle + \frac{1}{2\sigma^t}\|\bu^{t+1} - \bu^t\|^2$}
\State $\sigma^t = \underline \rho \sigma^t$, for any $\underline \rho <1$ and go back to Step 2.
\Else \State $\sigma^{t+1} = \overline \rho\sigma^t$, for any $\overline \rho >1$.
\EndIf
\end{algorithmic}
\caption{An iteration of Linearized Bregman iterations.}\label{alg:bregman_1v}
\end{algorithm}



\section{Numerical results}\label{sec:numerical_experiments}
We consider two sets of data, the first one, related to real measured data and a second one using synthetic (simulated) data. In both cases, we consider three energies labeled as $E_0$, $E_1$ and $E_2$, which are reconstructed separately. We analyze each energy channel independently as individual optimization problems. We compare the results from forward-backward splitting and linearized Bregman iterations and, highlight the main differences between TV and dTV regularizers. These algorithms were implemented using Python programming language and the Operator Discretization Library (ODL)~\cite{adler2017operator}. For each energy channel, we consider sinograms of size $60\times \text{ndet}$, \emph{i.e.} 60 projection angles and \emph{ndet} detectors. The angles are uniformly distributed in the interval $[0, 2\pi)$ and the reconstructed images $\bu$ are of size $512\times 512$. For measured data, $\text{ndet} = 552$ and, for synthetic, $\text{ndet} = 640$.

First, we detail how to choose the side information $\bv$ in our experiments considering the multi-spectral information in each energy channel.

\subsubsection{Choice of side information}\label{sec:choice_sinfo}

We propose to reconstruct a polyenergetic image $\bv\in \Rbb^N$ based on combining the data sets $\bb_k\in \Rbb^{M}$ for $k\in \{1, 2, 3\}$, \emph{i.e.}, we solve
\begin{equation}\label{eq:sideinfo_opt}
  \bv \in  \argmin_{\bx\in \Rbb^N} \left\{\frac 12\|A\bx-\tilde\bb\|^2_2 +\alpha \TV(\bx) + \iota_{[0,\infty)^N}(\bx)\right\}.
\end{equation}
where $\tilde\bb=\sum_{k=1}^3 \bb_k$. The regularization parameter $\alpha$ and more details related to this optimization problem will be specified during the numerical experiments for synthetic and real data. Solving~\eqref{eq:sideinfo_opt}, we get an image $\bv$ that despite of losing the spectral resolution, keeps structural information provided by all energy levels. Additionally, this image has higher signal-to-noise ratio and helps to improve the individual reconstructions $\bx_k$ as we show in our experiments.

We present the results using red, green and blue color maps for $E_0$, $E_1$ and $E_2$, respectively. We use the color grey to distinguish everything related to side information, making an analogy with the grayscale representation of an RGB image.

\subsection{Real data experiments}\label{sec:real_data}

Experimental data was gathered at the Department of Physics, University of Helsinki, using a cone-beam micro-CT scanner with an end-window tube and a tungsten target (GE Phoenix nanotom 180 NF). The phantom used in this experiment corresponds to the cross section of a small bird (at chest level) with different tissue types and fine bone structures. This phantom was previously introduced in~\cite{toivanen2019joint}, it was imaged using three different X-ray tube settings averaging four frames in each projection to increase signal-to-noise ratio. 
The three energy datasets were generated using 50, 80 and 120 kV, respectively. For $E_0$, an electric current of \SI{300}{\micro\ampere}, an exposure time of \SI{125}{\micro\second} and no filtration were used; for $E_1$, the experiment was setted up with \SI{180}{\micro\ampere}, \SI{125}{\micro\second} and \SI{1}{\milli\metre} of aluminium filtration. Similarly, for $E_2$, \SI{120}{\micro\ampere}, \SI{250}{\micro\second} and \SI{0.5}{\milli\metre} of copper filtration were configured. 2D sinograms were created using the central plane of the cone-beam projections, in which the geometry reduces to a fan-beam geometry. 
First, we discuss about the choice of reference images and side information.

\begin{figure}[htbp]
    \centering
    \includegraphics[width=0.98\linewidth]{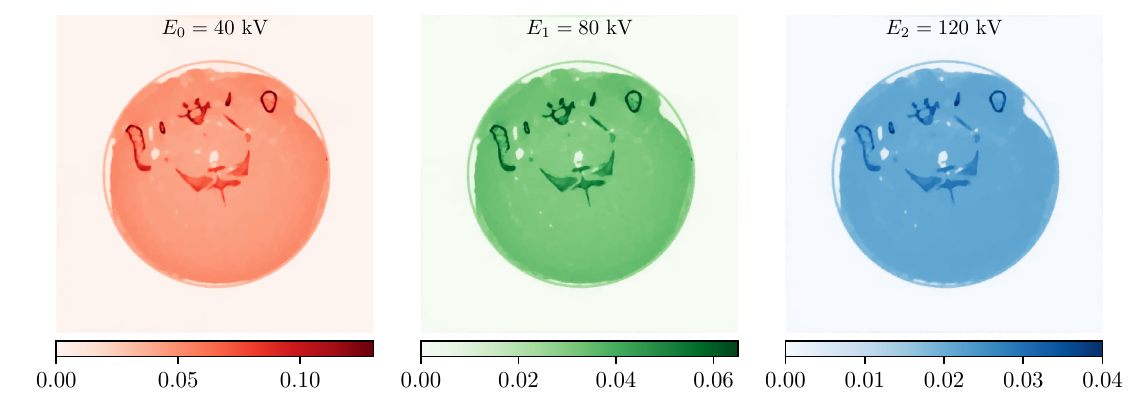}
    \caption{Reference images for real data solving the problem~\eqref{eq:variational_1v} with $\alpha=0.005$, $\alpha=0.002$ and $\alpha=0.002$ for $E_0$, $E_1$ and $E_2$, respectively. }\label{fig:real_groundtruths}
\end{figure}

\paragraph{Reference images:} Given that the object of study was the cross section of a bird and not precisely a human patient, a complete scan was carried out for each energy, that is, an experiment with high radiation dose and high exposure time to observe 720 rotation angles. This provided us a high resolution image to be used as ground truth. More precisely, the reference image for each energy uses 720 angles and 552 detectors and is computed via~\eqref{eq:variational_1v} using FBS and TV regularizer. The regularization parameter $\alpha$ is chosen to preserve low noise and high-resolution details in the reconstructions. The proposed references are in~\cref{fig:real_groundtruths}. Alternative references can be obtained using Filtered-backprojection algorithm, in our case, since the choice of optimal $\alpha$ (for FBS) and optimal iteration (Bregman iterations) is related to the PSNR values, the use of a FBP reconstruction as reference did not give us good results, on the contrary, images with few details or a lot of noise were obtained.

\begin{figure}[htbp]
  \centering
  \begin{tikzpicture}[spy using outlines={rectangle, frames,magnification=2.5, size=1cm, connect spies, every spy on node/.append style={ultra thick}}, x=\PicWidth, y=\PicWidth]
    \node {\pgfimage[width=0.87\linewidth]{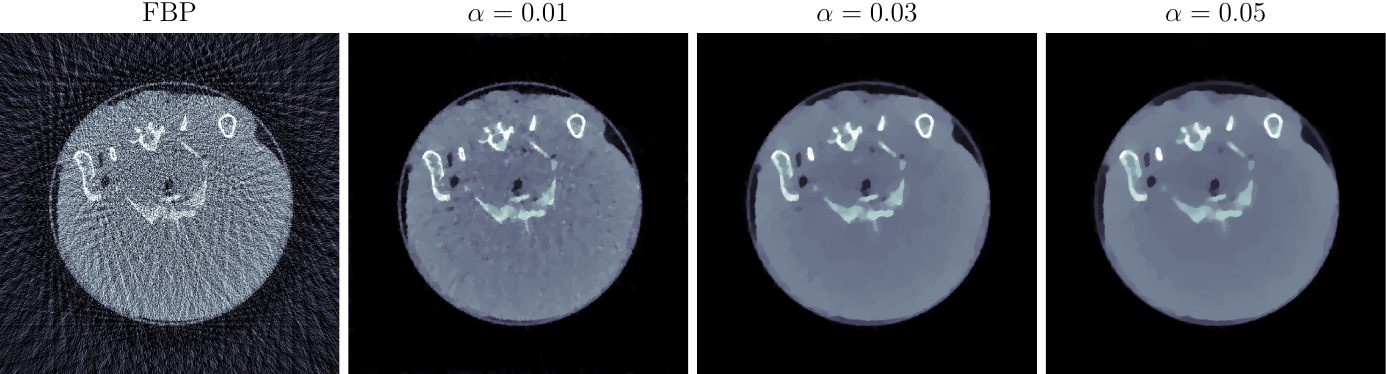}};
    \spy on (-5.8,0.1) in node [right] at (-4,1.5);
    \spy on (-2.4,0.1) in node [right] at (-0.6,1.5);
    \spy on (1,0.1) in node [right] at (2.8,1.5);
    \spy on (4.4,0.1) in node [right] at (6.2,1.5);
    \end{tikzpicture}\\
\caption{From left to right: FBP reconstruction of side information with 60 angles and 552 detectors, TV solutions of the problem~\eqref{eq:variational_1v} with $\alpha = 0.01$, $\alpha = 0.03$ and $\alpha = 0.05$, respectively. The minimization problem is solved with FBS algorithm in a space of size $512\times 512$.}\label{fig:real_sinfo}
\end{figure}

\paragraph*{Choice of side information: } 
We solve the problem~\eqref{eq:sideinfo_opt} for different values of $\alpha$, we compare the resulting reconstructions in~\cref{fig:real_sinfo}, also including a FBP-side information. 
We chose the reconstruction with $\alpha=0.03$, that keeps sharper boundaries and includes few artifacts during the reconstruction. We will compare the potential differences when a FBP- and a TV-side informations are used in synthetic data section (see~\cref{fig:sinfo_comparison}), we will mainly observe the improvements obtained by a regularized prior image.

\begin{figure}[htbp]
  \centering
  \includegraphics[width=\linewidth]{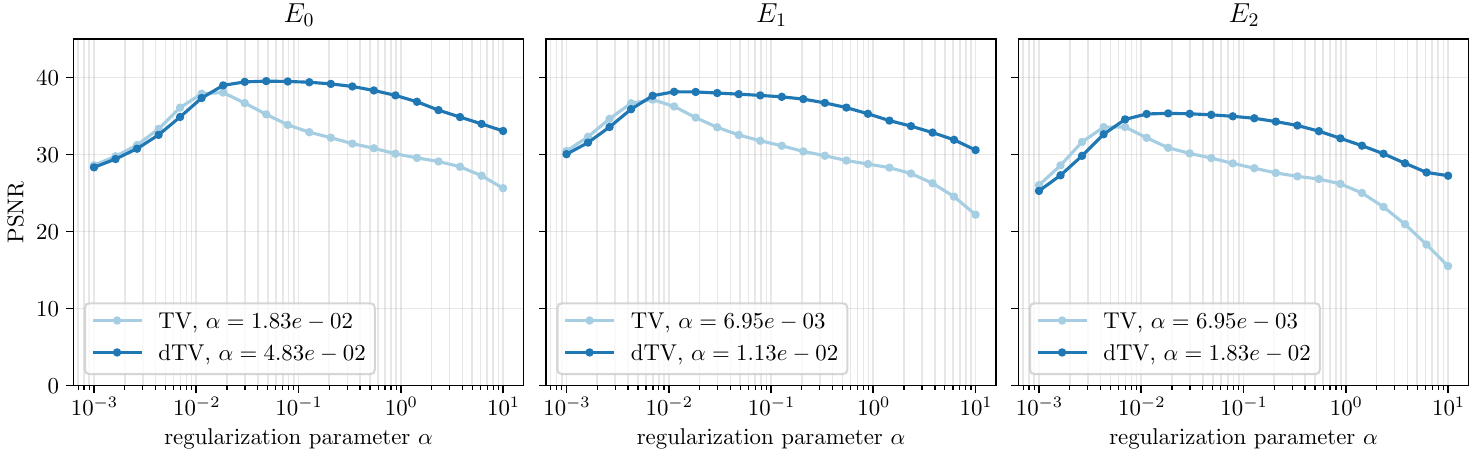}
  \caption{PSNR values for 20 different regularization parameters using 60 angles. FBS algorithm is used with TV and dTV regularizers and the optimal $\alpha$ corresponds to the highest value of PSNR.}\label{fig:bird_alphas_psnr}
  \end{figure}

\begin{figure}[htbp]
  \centering
  \begin{tabular}{@{}c@{ }c@{}}
    \begin{tikzpicture}[spy using outlines={rectangle,black,magnification=3,size=1.5cm, connect spies, every spy on node/.append style={ultra thick}}, x=\PicWidth, y=\PicWidth]
      \node {\pgfimage[width=0.87\linewidth]{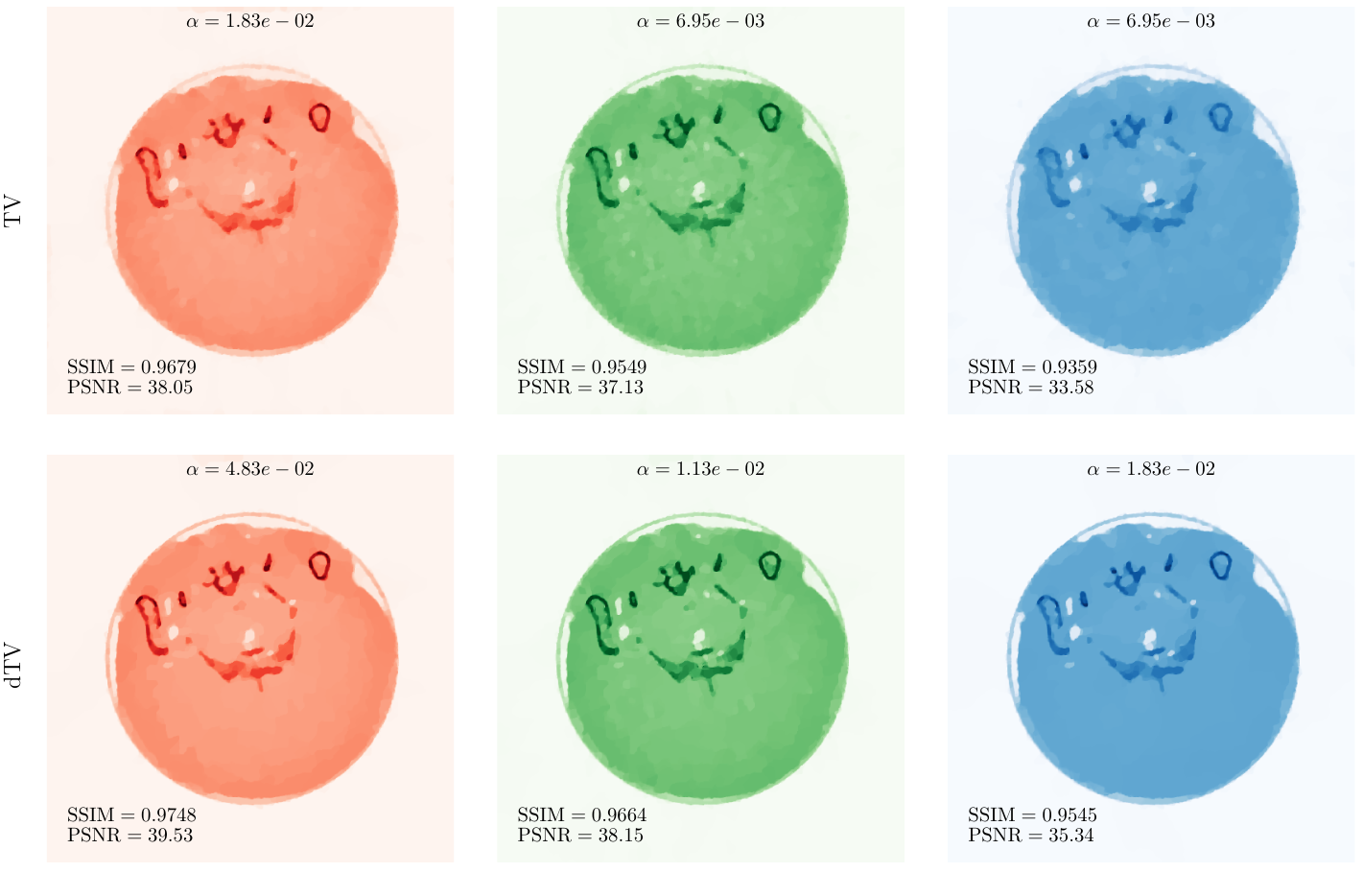}};
      \spy on (-5.1,2.8) in node [right] at (-3.3,4.2);
      \spy on (-0.635,2.8) in node [right] at (1.2,4.2);
      \spy on (3.83,2.8) in node [right] at (5.7,4.2);
      \spy on (-5.1,-1.63) in node [right] at (-3.3,-0.25);
      \spy on (-0.64,-1.63) in node [right] at (1.2,-0.25);
      \spy on (3.83,-1.63) in node [right] at (5.7,-0.25);
      \end{tikzpicture}\\
  \end{tabular}
  \caption{For the three energies, reconstructions using FBS algorithm with TV (upper row) and dTV (bottom row). For each setting, $\alpha$ is chosen to maximize PSNR.}\label{fig:bird_alphas_ssim}
  \end{figure}

\paragraph{FBS results:} We run the FBS iterations from~\cref{alg:fbs}, starting with $\bu^0=1\in \Rbb^N$ and $\sigma^0 = 1/\|A\|^2$, where $\|A\|$ is an estimated norm of the operator $A$. We set tolerance as $\tol=10^{-6}$ for all the experiments, so  the algorithm stops when $H(\bu^{t+1})-H(\bu^{t})\leq \tol\cdot H(\bu^{t+1})$. We choose $\eta=0.01\cdot \max_x |\nabla v(x)|$ for dTV definition~\eqref{eq:eta} as commonly done for this regularizer~\cite{ehrhardt2016multicontrast,bungert2018blind}. Since we are solving the problem~\eqref{eq:variational_1v}, the parameter $\alpha$ is chosen by running an arrange of values between $10^{-3}$ and $10^1$, the problem is solved for each value to compute the Peak Signal-to-Noise Ratio (PSNR)~\cite{hore2010image}, as in~\cref{fig:bird_alphas_psnr}. The optimal $\alpha$ corresponds to the highest PSNR. 
The results for all energies are shown in~\cref{fig:bird_alphas_ssim}, we have included the structural similarity measure (SSIM)~\cite{wang2004image} and PSNR measures implemented in ODL, to compare the quality of the reconstructions to the reference images in~\cref{fig:real_groundtruths}. 

\begin{figure}[htbp]
\begin{center}
\includegraphics[width=\linewidth]{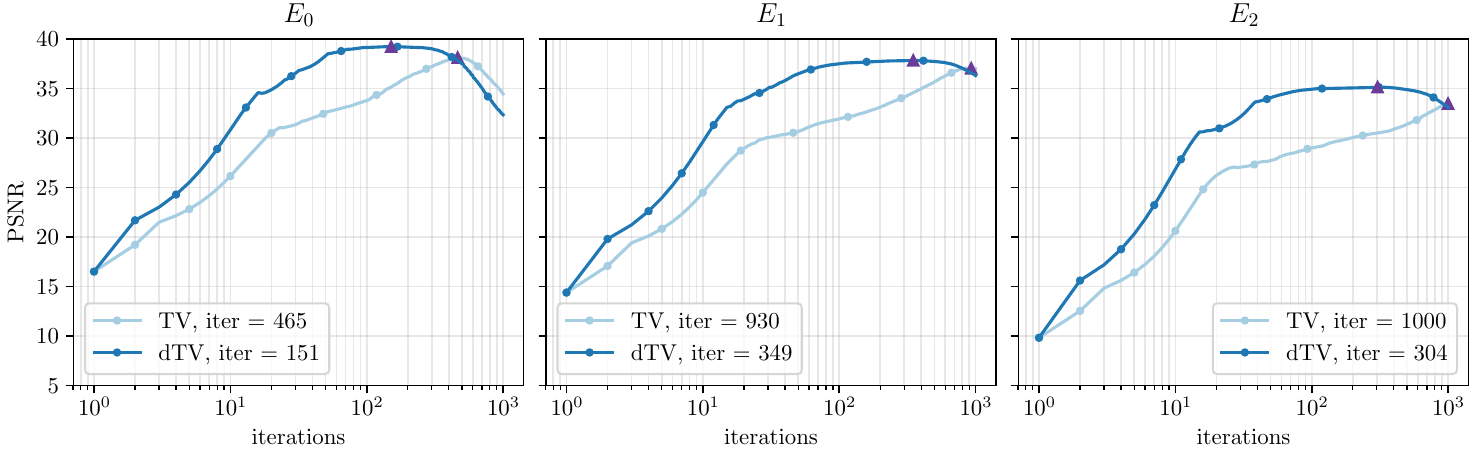}
\end{center}
\caption{For each energy level, the graphs show iterations against PSNR.}\label{fig:real_ssim_vs_iter}
\end{figure}

\begin{figure}[htbp]
  \centering
  \begin{tabular}{@{}c@{ }c@{}}
    \begin{tikzpicture}[spy using outlines={rectangle,black,magnification=3,size=1.5cm, connect spies, every spy on node/.append style={ultra thick}}, x=\PicWidth, y=\PicWidth]
      \node {\pgfimage[width=0.87\linewidth]{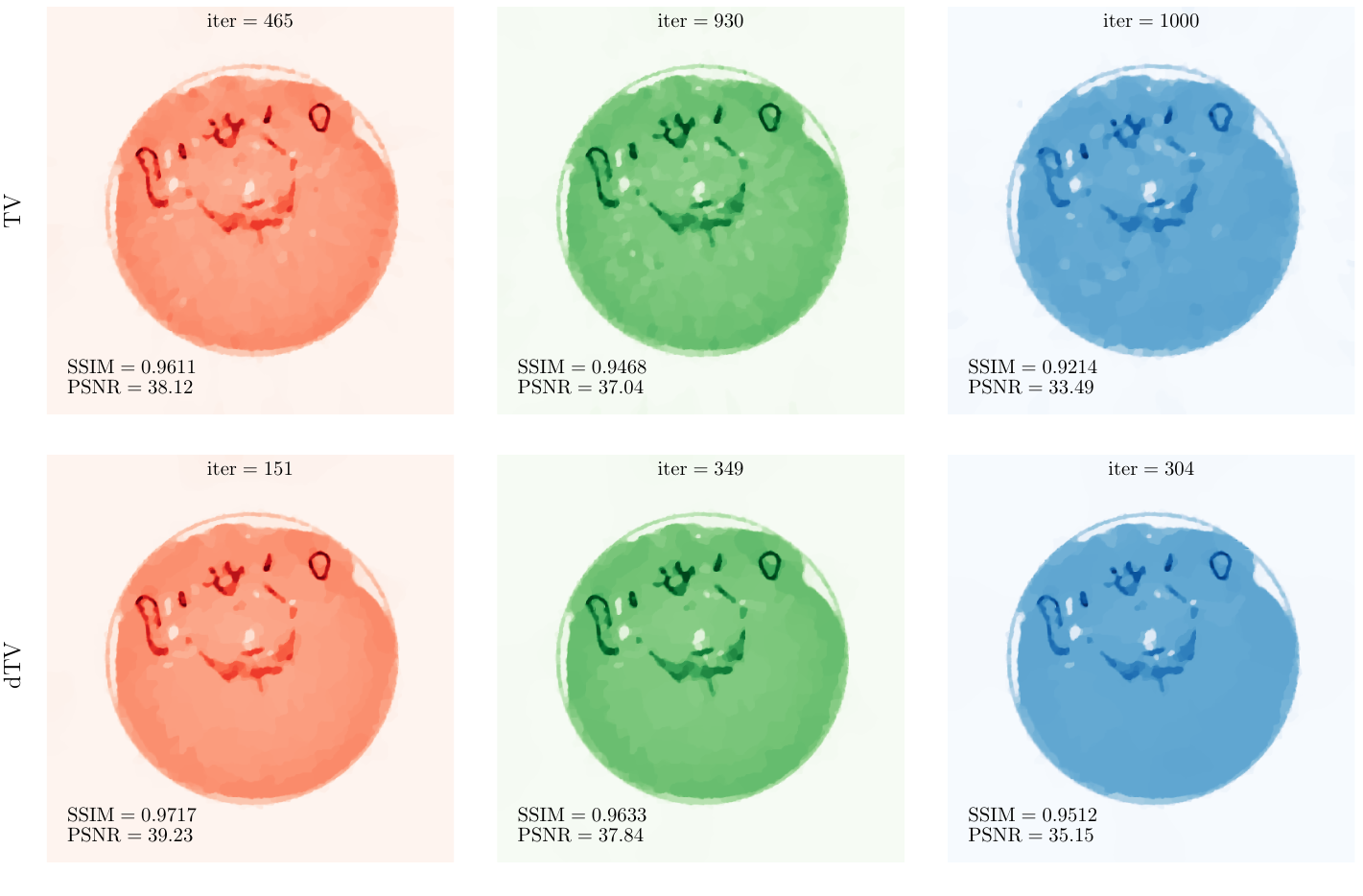}};
      \spy on (-5.1,2.7) in node [right] at (-3.3,4.2);
      \spy on (-0.635,2.7) in node [right] at (1.2,4.2);
      \spy on (3.83,2.7) in node [right] at (5.7,4.2);
      \spy on (-5.1,-1.73) in node [right] at (-3.3,-0.25);
      \spy on (-0.64,-1.73) in node [right] at (1.2,-0.25);
      \spy on (3.83,-1.73) in node [right] at (5.7,-0.25);
      \end{tikzpicture}\\
  \end{tabular} 
\caption{For the three energies, reconstructions using Bregman iterations with TV(upper row) and dTV (bottom row). Each image is labeled by the iteration that maximize PSNR.} \label{fig:real_all_energies_bregman}
\end{figure}

  \paragraph{Bregman results:} For \cref{alg:bregman_1v}, we start with $\sigma^0 = 1/\|A\|^2$ and $\bu^0 = \bq^0 = \bm 0\in \Rbb^N$. These choices guarantee that $\bq^0\in\partial G(\bu^0)$ for $G=G_{\TV}$ or $G=G_{\dTV}$. We run 1000 iterations in~\cref{alg:bregman_1v} using $\alpha=10$ for~\eqref{eq:G_TV} and \eqref{eq:G_dTV}. We observe from~\cref{fig:real_ssim_vs_iter} a common pattern along energies: the PSNR curve for TV is always below the curve associated to dTV, additionally, the number of iterations needed to maximize PSNR is always smaller for dTV than TV. The triangle markers refer to the best iterations in terms of PSNR. The images with highest PSNR are presented in~\cref{fig:real_all_energies_bregman}.

\begin{figure}[htbp]
    \centering
    \begin{tabular}{@{}c@{ }c@{}}
      \begin{tikzpicture}[spy using outlines={rectangle,black,magnification=2.5, size=1cm, connect spies, every spy on node/.append style={ultra thick}}, x=\PicWidth, y=\PicWidth]
        \node {\pgfimage[width=0.875\linewidth]{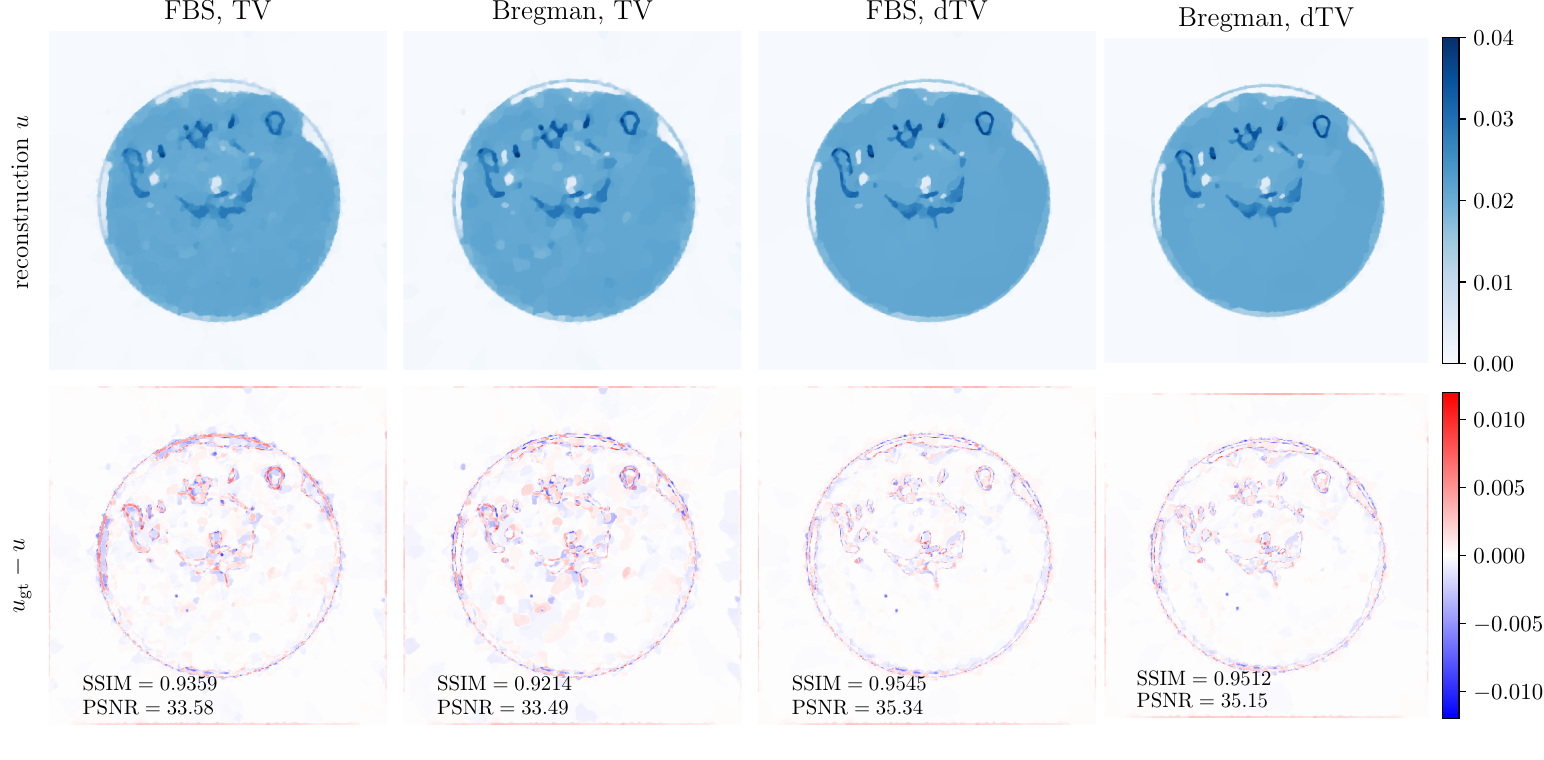}};
        \spy on (-5.45,1.8) in node [right] at (-4,2.5);
        \spy on (-2.34,1.8) in node [right] at (-1,2.5);
        \spy on (0.78,1.8) in node [right] at (2,2.5);
        \spy on (3.79,1.8) in node [right] at (5,2.5);
        \spy on (-5.45,-1.35) in node [right] at (-4,-0.3);
        \spy on (-2.34,-1.35) in node [right] at (-1.2,-0.3);
        \spy on (0.78,-1.35) in node [right] at (2,-0.3);
        \spy on (3.79,-1.35) in node [right] at (5,-0.3);
        \end{tikzpicture}\\
    \end{tabular} 
    \caption{Top: reconstructions for $E_2$ using ``optimal'' regularization for both FBS and Bregman iterations. Bottom: difference between reconstruction and reference image (see~\cref{fig:real_groundtruths}).}
    \label{fig:real_fbs_bregman_best}
\end{figure}

\paragraph{Comparison between FBS and Bregman iterations:}
After calibration of the regularization parameter $\alpha$ for FBS and iteration number for Bregman iterations, we compare the two algorithms for $E_2$ in~\cref{fig:real_fbs_bregman_best}. We observe that both algorithms give similar values for PSNR and SSIM but consistently dTV outperforms TV.

\subsection{Synthetic data experiments}
\begin{figure}[htbp]
  \begin{center}
  \includegraphics[width=0.55\linewidth]{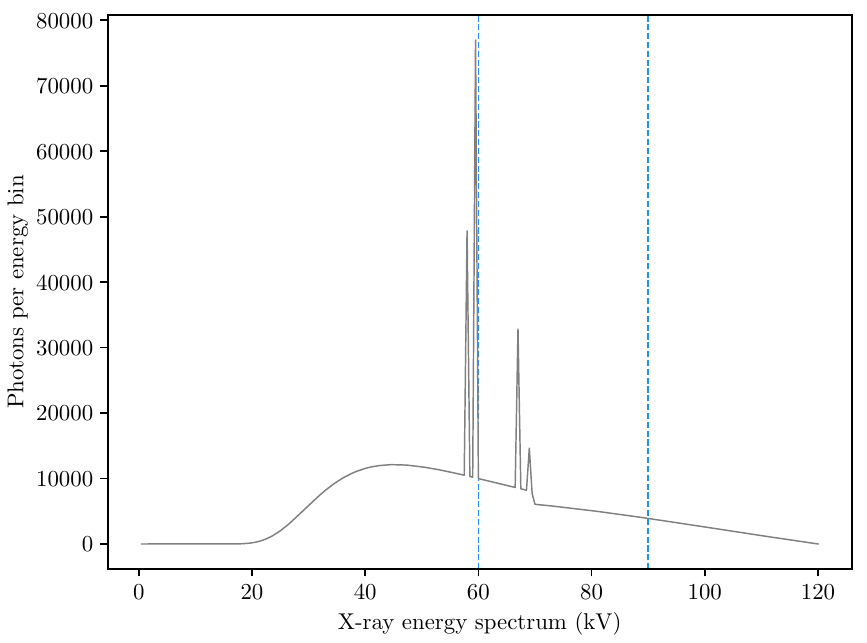}
  \caption{The energy spectrum with the three energy bins divided at \SI{60}{\kilo\volt} and \SI{90}{ \kilo\volt}.}\label{fig:spectrum_material}
  \end{center}
\end{figure}

\begin{figure}[htbp]
  \centering
    \begin{tikzpicture}[spy using outlines={rectangle,black,magnification=3,size=1.5cm, connect spies, every spy on node/.append style={ultra thick}}, x=\PicWidth, y=\PicWidth]
      \node {\pgfimage[width=0.875\linewidth]{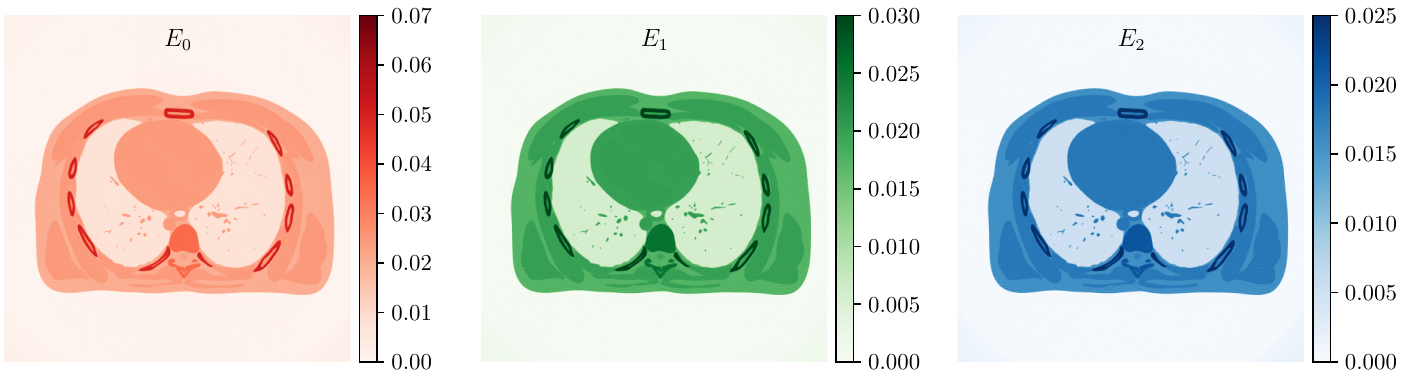}};
      \spy on (-4.4,0.3) in node [left] at (-5.5,1.3);
      \spy on (-0.4,-0.7) in node [left] at (-0.8,1.3);
      \spy on (4.6,-0.3) in node [left] at (3.8, 1.3);
      \end{tikzpicture}
  \caption{\emph{From left to right:} The energy spectrum with the three energy bins divided at \SI{60}{\kilo\volt} and \SI{90}{ \kilo\volt} and, reference images for $E_0$, $E_1$ and $E_2$. Each image presents a different zoomed zone as a reference for the oncoming experiments.} \label{fig:synt_gt}
  \end{figure}

For the simulated data, we used an anthropomorphic chest phantom generated using the XCAT software~\cite{segars2008realistic}. We simulated measurements using the XCAT CT Simulator, using an X-ray spectrum included with the software. The simulated spectrum was created using SRS-78 Spectrum Processor program, and it modelled a tungsten tube X-ray source with a 120 kV acceleration voltage and 5 mm of aluminum filtration. For the multi-energy measurement simulations, the spectrum was divided into three energy bins: 0-\SI{60}{\kilo\volt}, 60-\SI{90}{ \kilo\volt}, and 90-\SI{120}{\kilo\volt}. For testing the algorithms, we created 720 projections at 0.5 degree intervals, with simulated photon noise. For  creating the reference/ground truth, we simulated 1440 projections at 0.25 degree intervals, with no noise. In~\cref{fig:spectrum_material}, we present the energy spectrum for our data and in~\cref{fig:synt_gt} the reference images.

\paragraph*{Choice of side information: } 
As in real case, we solve the problem~\eqref{eq:sideinfo_opt} using  $\alpha = 0.08$, $\alpha = 0.1$ and $\alpha = 0.5$. We have chosen the image with $\alpha=0.1$, which gave us the best results. We discuss the accuracy of side information choice in~\cref{fig:sinfo_comparison}.

\begin{figure}[htbp]
\centering
\includegraphics[width=\linewidth]{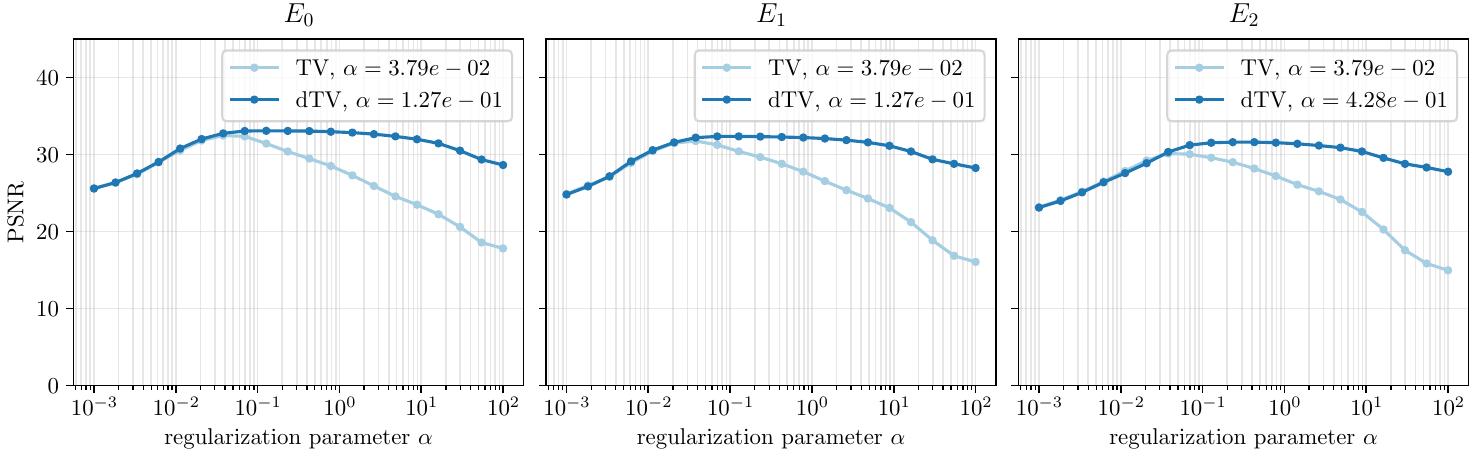}
\caption{PSNR values for 20 different regularization parameters using 60 angles. FBS algorithm is used with TV and dTV regularizers and the optimal $\alpha$ corresponds to the highest value of PSNR.}\label{fig:materials_alphas_ssim}
\end{figure}
\begin{figure}[htbp]
  \centering
  \begin{tabular}{@{}c@{ }c@{}}
    \begin{tikzpicture}[spy using outlines={rectangle,black,magnification=3,size=1.5cm, connect spies, every spy on node/.append style={ultra thick}}, x=\PicWidth, y=\PicWidth]
      \node {\pgfimage[width=0.87\linewidth]{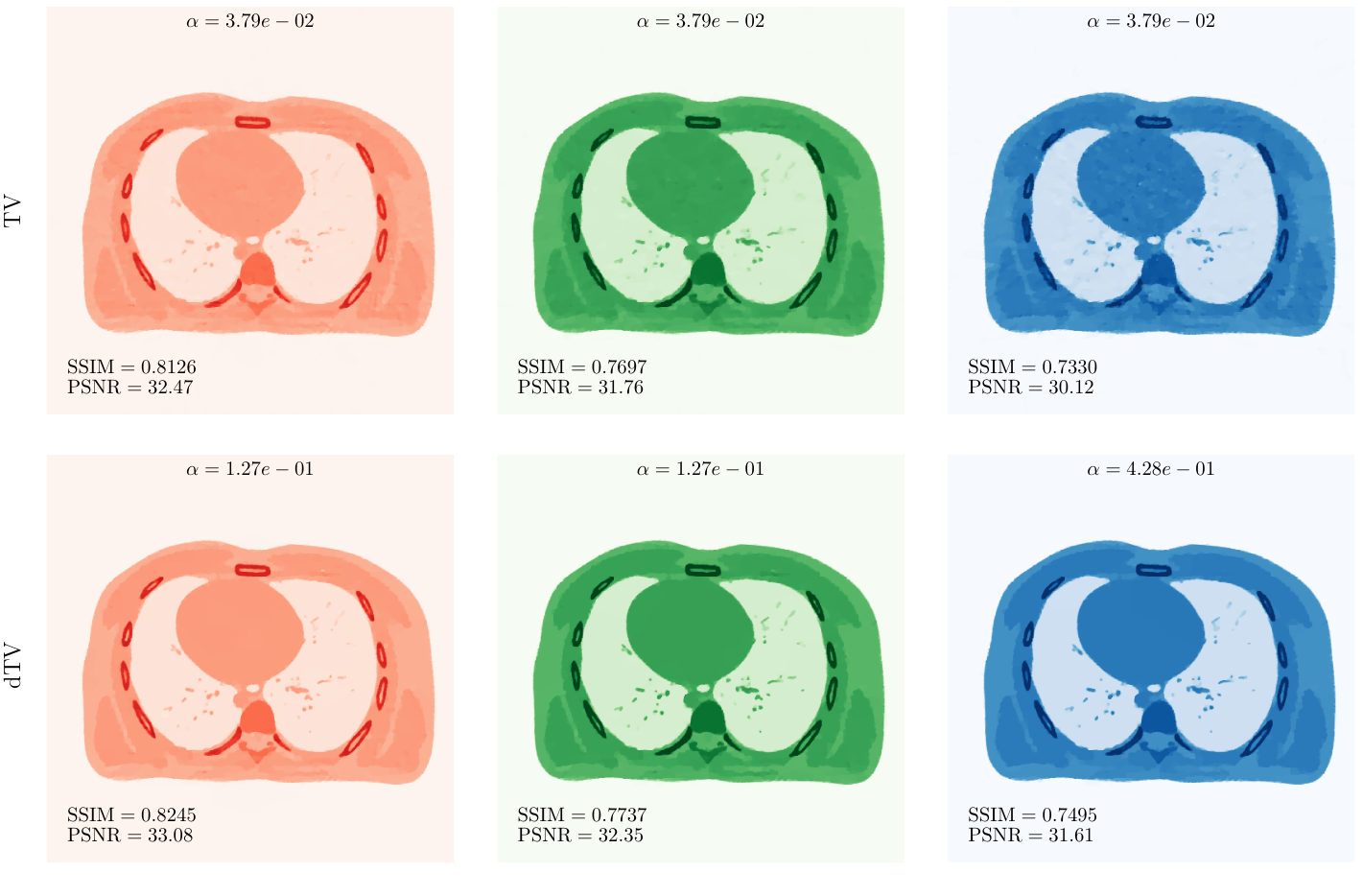}};
      \spy on (-3.6,2.7) in node [right] at (-3.3,4.2);
      \spy on (0.865,2.7) in node [right] at (1.2,4.2);
      \spy on (5.33,2.7) in node [right] at (5.7,4.2);
      \spy on (-3.6,-1.73) in node [right] at (-3.3,-0.25);
      \spy on (0.865,-1.73) in node [right] at (1.2,-0.25);
      \spy on (5.33,-1.73) in node [right] at (5.7,-0.25);
      \end{tikzpicture}
  \end{tabular} 
\caption{For the three energies, reconstructions using Bregman iterations with TV (upper row) and dTV (bottom row). Each image is labeled by the iteration that maximize PSNR.
} \label{fig:synthetic_all_energies_bregman}
\end{figure}

\paragraph{FBS results:} We initialized $\bu^0$ and $\sigma^0$ as in real experiments. In~\cref{fig:materials_alphas_ssim}, we present the PSNR values obtained for 20 different values of $\alpha$ from $10^{-3}$ to $10^2$ for TV and dTV regularizers in each energy channel. The corresponding reconstructions are in \cref{fig:synthetic_all_energies_bregman}.
We included close-ups for easier comparison of the reconstructions. Here, we observe that using dTV yields higher values of measures and better reconstructions reducing the artifacts in the zoomed area at the top.

\begin{figure}[htbp]
  \centering
  \includegraphics[width = 0.48\linewidth]{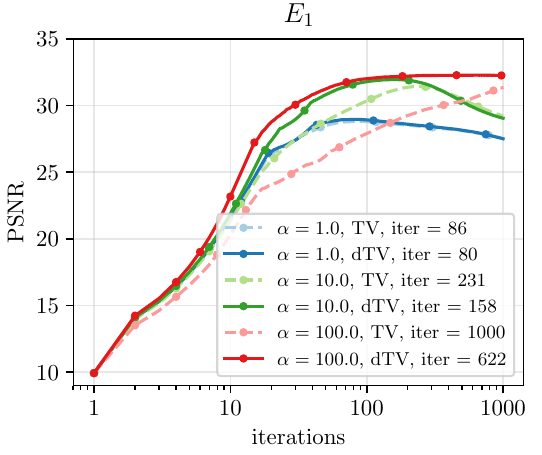}
  \includegraphics[width = 0.48\linewidth]{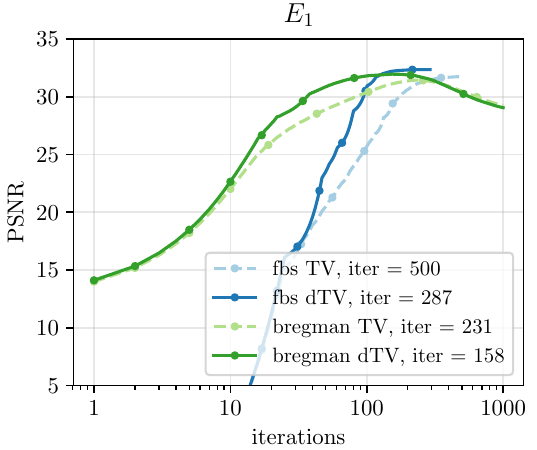}
  \caption{\emph{Left:} PSNR values for Bregman iterations using $\alpha = 1$, $\alpha =10$, and $\alpha=100$. The results are shown for both regularizers and the optimal iteration is included based on the highest PSNR value. \emph{Right: } PSNR along iterations for FBS and Bregman iterations using both TV and dTV.}\label{fig:E1_alphas_bregman}
\end{figure}

\begin{figure}[htbp]
  \centering
    \begin{tikzpicture}[spy using outlines={rectangle,black,magnification=3,size=1.8cm, connect spies, every spy on node/.append style={ultra thick}}, x=\PicWidth, y=\PicWidth]
      \node {\pgfimage[width=0.3915\linewidth]{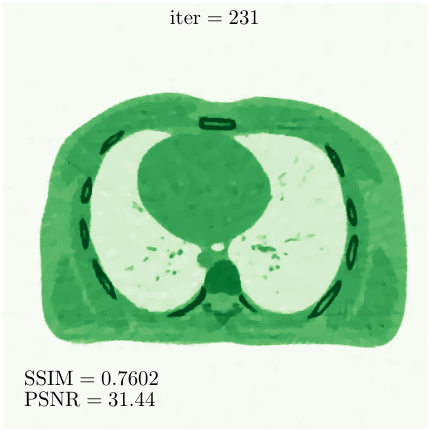}};
      \spy on (0.25,-1.4) in node [right] at (1.5,2.5);
      \spy on (-0.78,-0.7) in node [left] at (-1.5,2.5);
      \end{tikzpicture}
      \begin{tikzpicture}[spy using outlines={rectangle,black,magnification=3,size=1.8cm, connect spies, every spy on node/.append style={ultra thick}}, x=\PicWidth, y=\PicWidth]
        \node {\pgfimage[width=0.3915\linewidth]{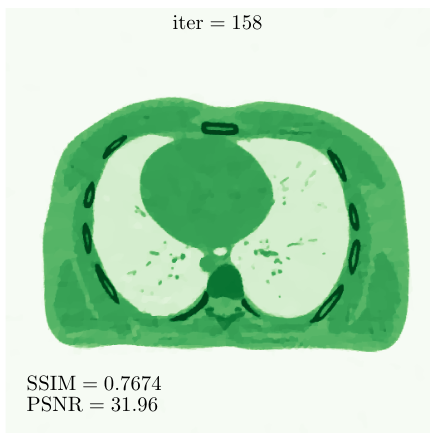}};
        \spy on (0.25,-1.4) in node [right] at (1.5,2.5);
        \spy on (-0.78,-0.7) in node [left] at (-1.5,2.5);
        \end{tikzpicture}\\
\caption{Bregman iterations for TV (left) and dTV (right) for $E_1$ using $\alpha=10$ based on the left hand of~\cref{fig:E1_alphas_bregman}.}\label{fig:E1_bregman_TV_dTV}
\end{figure}

\paragraph{Bregman results:} Now, we observe the results for~\cref{alg:bregman_1v}. We initialize $\sigma^0 = 1/\|A\|^2$ and $\bu^0 = \bq^0 = \bm 0\in \Rbb^N$. We run the algorithm for 1000 iterations using three different regularization parameters,  $\alpha = 1$, $\alpha =10$, and $\alpha=100$. As before, the ``optimal'' iteration number is chosen to maximize PSNR. In~\cref{fig:E1_alphas_bregman} (left), using $E_1$ as an example, we present the PSNR values along iterations for TV and dTV when these three values of $\alpha$ are considered in the minimization functional. We observe the influence of $\alpha$ in terms of number of iterations and quality of the reconstruction. While choosing a small value of $\alpha$ allows us to reach a maximum value of PSNR in a small number of iterations, the PSNR values are smaller compared to a higher $\alpha$ value. The highest values of PSNR occur when $\alpha=100$ but TV reconstructions need more than 1000 of iterations to find a maximum PSNR. In~\cref{fig:E1_bregman_TV_dTV}, we present the resulting reconstructions for TV and dTV in $E_1$. We have chosen $\alpha=10$ since the PSNR values in this case are close to the ones obtained with $\alpha=100$ and the number of iterations are less than 250.

The FBS and Bregman reconstructions for $E_1$ in~\cref{fig:synthetic_all_energies_bregman} and \cref{fig:E1_bregman_TV_dTV}  reach similar results, however, we note that fewer iterations were needed for Bregman iterations compared to FBS as shown in~\cref{fig:E1_alphas_bregman} (right). This means that linearized Bregman iterations converge faster to desired solution than forward-backward splitting.

\subsubsection{Influence of side information}
\begin{figure}[htbp]
  \centering
  \begin{tabular}{@{}c@{ }c@{}}
    \begin{tikzpicture}[spy using outlines={rectangle, frames, magnification=3,size=2cm, connect spies, every spy on node/.append style={ultra thick}}, x=\PicWidth, y=\PicWidth]
      \node {\pgfimage[width=0.88\linewidth]{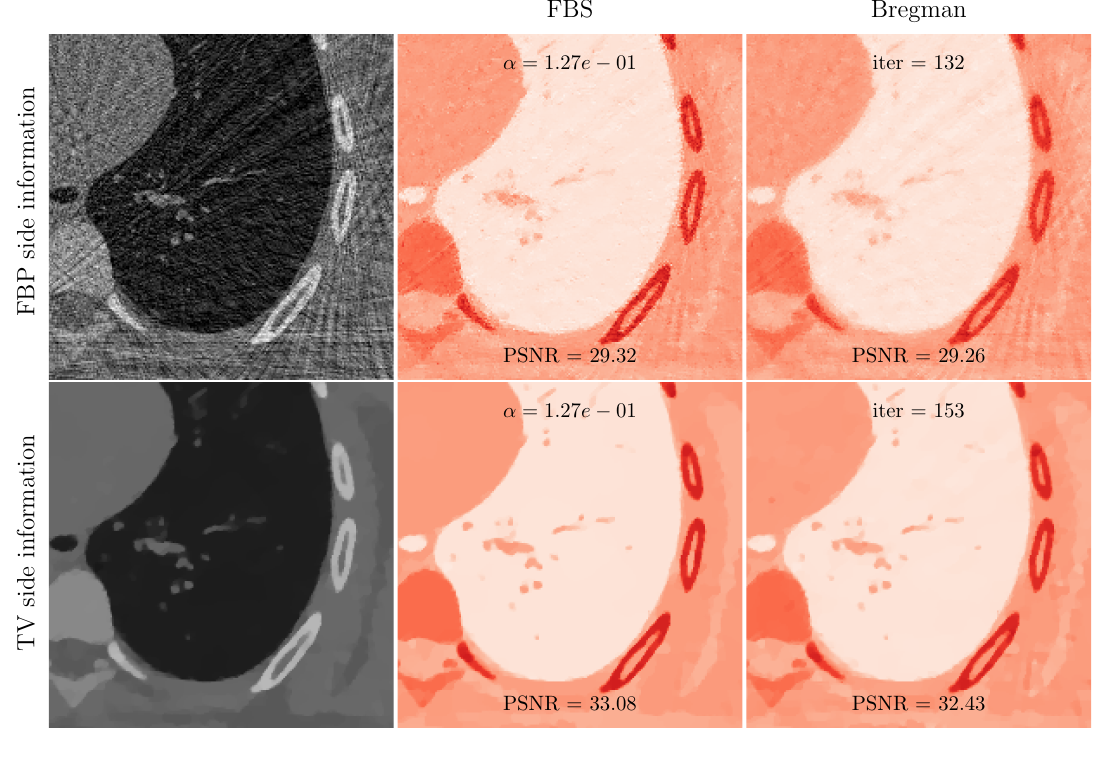}};
      \spy on (-4.5,2) in node [right] at (-3.55,3.7);
      \spy on (-0.23,2) in node [right] at (1.05,3.7);
      \spy on (4.0,2) in node [right] at (5.4,3.7);
      \spy on (-4.5,-2.3) in node [right] at (-3.55,-0.35);
      \spy on (-0.23,-2.3) in node [right] at (1.05,-0.35);
      \spy on (4.0,-2.3) in node [right] at (5.4,-0.35);
      \end{tikzpicture}
  \end{tabular}  
    \caption{Comparison between dTV reconstructions using FBS and Bregman iterations with two different side informations, one TV-regularized with $\alpha=100$ and a second, non-regularized side information with FBP.}
    \label{fig:sinfo_comparison}
\end{figure}

In this experiment, we compare the accuracy of the reconstructions depending on the choice of side information. For this, we consider a FBP- an TV- side information reconstructions. The second one was obtained by solving~\eqref{eq:sideinfo_opt} with $\alpha = 0.1$. These two choices, shown in~\cref{fig:sinfo_comparison} for XCAT phantom, give us two types of side information, one image with more structures and artifacts (upper row) and another image with smoother shapes but without artifacts (bottom row). We compare the best reconstructions using dTV with FBS algorithm considering the highest values of PSNR. For the first side information, the artifacts remain in the reconstructed image and poor quality values of PSNR are reached using larger number of Bregman iterations compared to the second case.

\section{Conclusion}
We have analyzed synergistic reconstruction for multi-spectral CT reconstruction when a limited set of angles is observed. The proposed approach is based on combining information from all available energy channels into a polyenergetic image. This image is then included into the directional total variation regularizer for use in variational or iterative regularization.

We observed that the synergistic approach based on directional total variation is always superior to separate reconstruction using just total variation for both variational and iterative regularization. In addition, we consistently saw that linearized Bregman iterations converge faster to a desired solution than forward-backward splitting.

Additionally, reconstructing images channel-by-channel due to dTV with channel-independent side information allows us to access the computational advantages of parallel programming, differentiating our proposal from methods that solve the multi-spectral problem by a single large optimization problem.

The observation that synergistic reconstruction can be faster than separate reconstruction is novel and interesting and future work will be directed to fully understand this phenomenon.

\subsubsection*{Appendix: Comparison with TVN}
\begin{figure}[htbp]
    \centering
    \begin{tikzpicture}[spy using outlines={rectangle, black, magnification=2, size=1.5cm, height=3.5cm, connect spies, every spy on node/.append style={ultra thick}}, x=\PicWidth, y=\PicWidth]
      \node {\pgfimage[width=0.8265\linewidth]{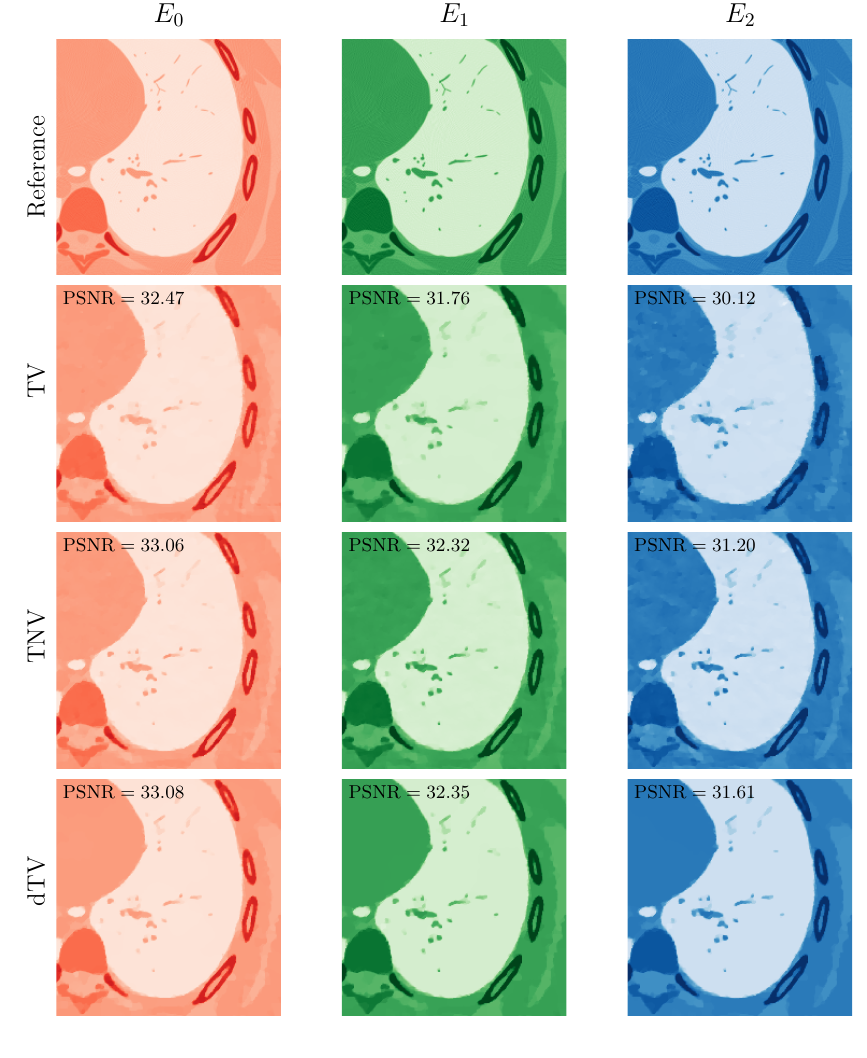}};
      \spy on (-4.9, 5) in node [right] at (-2.9,5.6);
      \spy on (-0.65,5) in node [right] at (1.35,5.6);
      \spy on (3.6,5) in node [right] at (5.6,5.6);

      \spy on (-4.9,1.35) in node [right] at (-2.9,1.95);
      \spy on (-0.65,1.35) in node [right] at (1.35,1.95);
      \spy on (3.6,1.35) in node [right] at (5.6,1.95);
      
      \spy on (-4.9,-2.3) in node [right] at (-2.9,-1.7);
      \spy on (-0.65,-2.3) in node [right] at (1.35,-1.7);
      \spy on (3.6,-2.3) in node [right] at (5.6,-1.7);

      \spy on (-4.9,-5.95) in node [right] at (-2.9,-5.35);
      \spy on (-0.65,-5.95) in node [right] at (1.35,-5.35);
      \spy on (3.6,-5.95) in node [right] at (5.6,-5.35);
      \end{tikzpicture}
    \caption{Comparison between TV, TVN and dTV regularizers. The PSNR values are included in each reconstruction. The noiseless images (references) for each energy are included in the top row. }
    \label{fig:comparison_dTV}
\end{figure}

Comparing the results with other methodologies used for multi-spectral CT is pertinent, for example with total nuclear variation (TNV) described in~\cite{rigie2017assessment} on realistic dual-energy CT data.  To do this, we solved the problem:
\[
\bu^\ast \in  \argmin_{\bu\in (\Rbb^N)^3} \left\{\sum_{i=0}^2\|R\bu_i-\bb_i\|^2_2 +\alpha \text{TNV}(\bu) + \iota_{[0,\infty)^N}(\bu)\right\}.
\]
where $\bu = (\bu_0, \bu_1, \bu_2)$ correspond to jointly find the reconstruction for the three images for $E_0$, $E_1$ and $E_2$. The results obtained are presented and compared with TV and dTV in Figure~\ref{fig:comparison_dTV}. The optimization problem was solved using Primal-dual hybrid gradient (pdhg) algorithm described in~\cite{chambolle2011first}. More details related to TNV with pdhg can be found in~\cite{Ehrhardt2020chapter}.

Comparing PSNR values for TNV and dTV, we can observe that the differences are quite negligible. In the case of dTV, sharper edges can be observed compared to TV and TNV. We could conclude that TNV and dTV methods give similar results, however, there are important differences in terms of computational cost. In TVN a large problem needs to be solved, if we have to reconstruct $N$ energy channels all have to be solve in the same computing unit, while using dTV in a channel-by-channel reconstruction allows us to solve these $N$ problems in $N$ different computers (nodes) through a  parallelization process.

\vskip6pt

\enlargethispage{20pt}


\noindent\textbf{Authors’ Contributions: }{EC carried out the numerical computations and drafted the manuscript. MJE designed the project, advised EC and supported writing the manuscript. SS and AM acquiered and processed real data. All authors read and approved the manuscript.}

\noindent\textbf{Funding: }{MJE acknowledges support from the EPSRC (EP/S026045/1, EP/T026693/1), the Faraday Institution (EP/T007745/1) and the Leverhulme Trust (ECF-2019-478). EC acknowledges support from the CMM ANID PIA AFB170001 and Beca
Doctorado Nacional Conicyt.}



\end{document}